\documentclass{article}

\usepackage[english]{babel}

\usepackage[utf8]{inputenc}
\usepackage[T1]{fontenc}

\usepackage{amssymb}
\usepackage{stmaryrd}
\usepackage{amsmath,amsfonts}
\usepackage[tikz]{bclogo}
\usepackage[top=4.34cm]{geometry}

\usepackage{lmodern}
\usepackage{textcomp}
\usepackage{mathtools, bm}
\usepackage{amssymb, bm}
\usepackage[unq]{unique}
\usepackage{enumerate}

\usepackage[breaklinks]{hyperref}
\usepackage[numbers]{natbib}    

\usepackage{amsthm}
\usepackage{comment}
\usepackage{todonotes}

\newtheorem{theorem}{Theorem}[section]
\newtheorem{lemma}{Lemma}[section]
\newtheorem{Obs}[theorem]{Observation}
\newtheorem{conjecture}{Conjecture}
\newtheorem{Prop}[theorem]{Proposition}

\theoremstyle{definition}

\title{A new lower bound on the pebbling number of the grid}
\author{Jan Petr\footnote{\href{mailto:jp895@cam.ac.uk}{jp895@cam.ac.uk}}, Julien Portier\footnote{\href{mailto:jp899@cam.ac.uk}{jp899@cam.ac.uk}}, Szymon Stolarczyk\footnote{\href{mailto:ss394470@students.mimuw.edu.pl}{ss394470@students.mimuw.edu.pl}}}

\begin{document}

\maketitle

\begin{abstract}
    A pebbling move on a graph consists of removing $2$ pebbles from a vertex and adding $1$ pebble to one of the neighbouring vertices. A vertex is called reachable if we can put $1$ pebble on it after a sequence of moves. The optimal pebbling number of a graph is the minimum number $m$ such that there exists a distribution of $m$ pebbles so that each vertex is reachable. For the case of a square grid $n \times m$, Gy\H{o}ri, Katona and Papp recently showed that its optimal pebbling number is at least $\frac{2}{13}nm \approx 0.1538nm$ and at most $\frac{2}{7}nm +O(n+m) \approx 0.2857nm$. We improve the lower bound to $\frac{5092}{28593}nm +O(m+n) \approx 0.1781nm$.
\end{abstract}

\section{Introduction}

Let $G$ be a graph. A \emph{pebbling distribution} is a function $P$ from $V(G)$ to the set of non-negative integers. We say that a vertex \emph{$x \in V(G)$ has $k$ pebbles on it} if $P(x)=k$. The \emph{total number of pebbles} is defined as $|P|=\sum_{x \in V(G)} P(x)$. A \emph{pebbling move} consists of removing $2$ pebbles from a vertex and adding $1$ pebble to one of the neighbouring vertices. We say that a vertex is \emph{reachable} if we can put at least $1$ pebble on it after a sequence of moves. More generally, for any $k \in \mathbb{N}$ we say that a vertex is \emph{$k$-reachable} if we can put at least $k$ pebbles on it after a sequence of moves. We say that a vertex is \emph{exactly $k$-reachable} if it is $k$-reachable but not $(k+1)$-reachable. We say that a pebbling distribution is \emph{solvable} if every vertex is reachable. The \emph{optimal pebbling number $\pi(G)$ of $G$} is the minimum number $m$ such that there exists a solvable distribution of $m$ pebbles on the graph $G$. 

Pebbling was introduced by Chung \cite{Chung} in 1989; since then the optimal pebbling number has been studied for various classes of graphs, such as paths and cycles (see, for example, \cite{Bundeetal}), $m$-ary trees \cite{mtrees} or hypercubes \cite{LowerBoundQn}, \cite{UpperBoundGrid}. As shown by Milans and Clark, the problem of determining $\pi(G)$ for a general graph $G$ is NP-complete \cite{NPC}. \\

In this paper, we will focus mainly on the case where $G$ is a $m \times n$ grid $\Lambda_{m,n}$: the vertices are the $nm$ squares of the grid, and two vertices share an edge if and only if their respective squares are adjacent. The best known (and conjectured optimal) upper bound is $\pi(\Lambda_{m,n}) \leq \frac{2}{7}nm+O(m+n) \approx 0.2857nm $, see a paper by Gy\H{o}ri, Katona and Papp \cite{UpperBoundGrid} for an explicit construction. In \cite{LowerBoundGrid}, the same authors proved that $\pi(\Lambda_{m,n}) \geq \frac{2}{13}nm \approx 0.1538nm$. We improve this result to $\pi(\Lambda_{m,n}) \geq \frac{5092}{28593}nm+O(m+n) \approx 0.1781nm$. \\

The core method used in our new lower bound has been used before, for example in \cite{LowerBoundGrid}. Given a graph $G$ and a pebbling distribution $P$, we define for each vertex $y$ the \emph{contribution function of $y$} as $v_y:V(G) \rightarrow \mathbb{R}$ given by $v_y(x)=P(y)2^{-d(x,y)}$ where $d$ is the graph distance in $G$. We also define the \emph{effect of a pebble placed on $x$} as $\mathrm{ef}(x)=\sum_{y \in V(G)} 2^{-d(x,y)}$. We define the \emph{value of $x$} as $v(x)=\sum_{y \in V(G)} v_y(x)$. The basic observation is that if $P$ is solvable, then the value of each vertex in $G$ is at least $1$. This already enables us to get a lower bound on $|P|$, using $\sum_{x \in V(G)}|P(x)|\mathrm{ef}(x) = \sum_{x \in V(G)} v(x) \geq |V(G)|$. Note that for $x \in V(\Lambda_{m,n})$ we have $\mathrm{ef}(x) \leq 1+\sum_{k=1}^{\infty} 4k\frac{1}{2^k}=9$. This bound already gives $\pi(\Lambda_{m,n}) \geq \frac{1}{9} mn$. In order to improve this, we will give a better lower bound on $\sum_{x \in V(G)} v(x)$. \\

To the best of our knowledge, our first improvement is a new concept in the study of pebbling. Given a solvable pebbling distribution, we partition the vertices into what we will call \emph{regions} in a way that all vertices inside a region can be reached using pebbles from within the region and no pebble can ever leave a region: given a graph $G$ and a starting distribution $P$, a \emph{region of reachability under $P$} (\emph{region} for short) is the set of vertices of a maximal connected subgraph of $G$ on $2$-reachable vertices together with their neighbours.

\begin{figure}[htbp]\centering
    			\includegraphics[height=7cm]{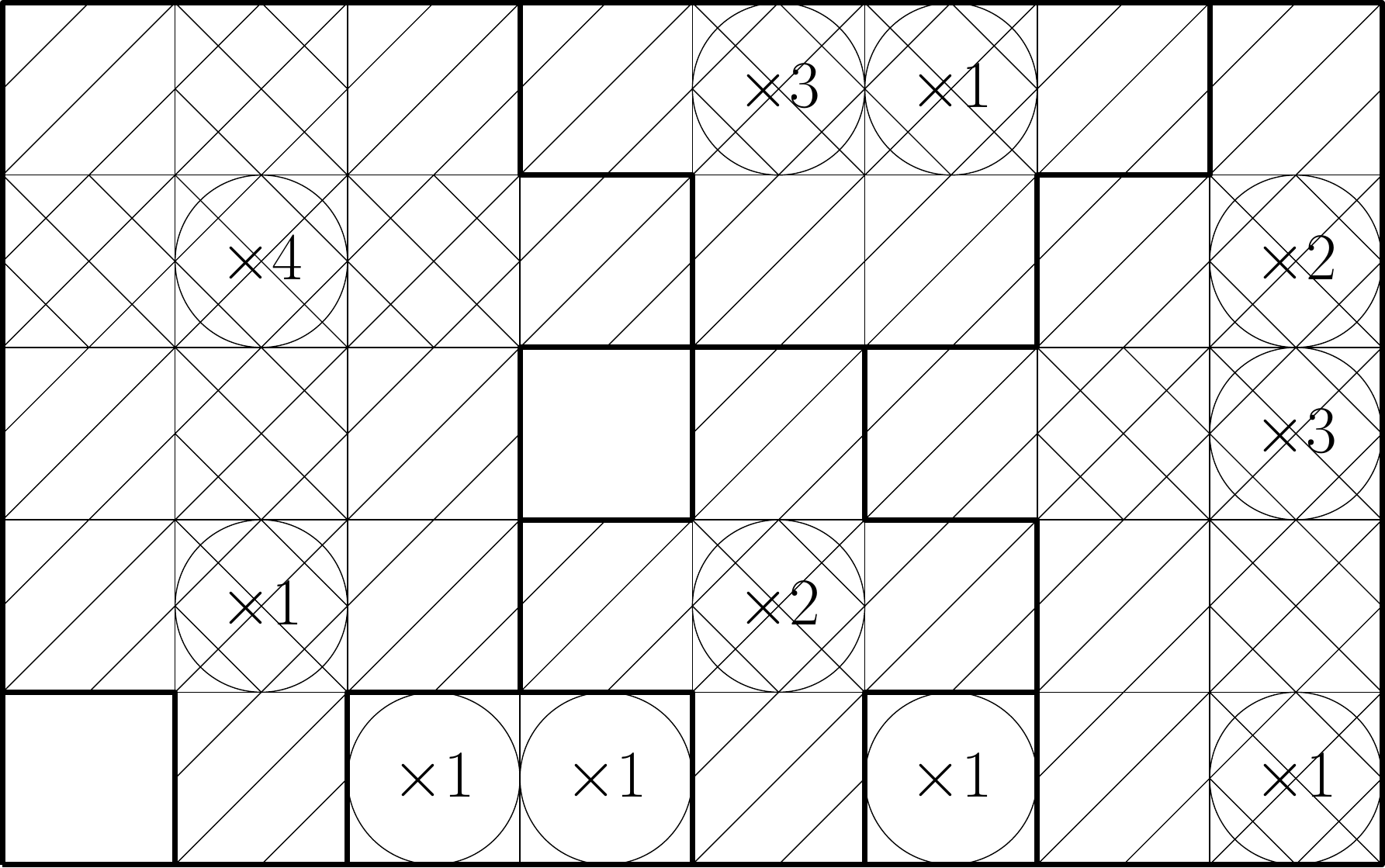}
    			\caption{A starting distribution $P$ on $\Lambda_{5,8}$ and the four regions of reachability under $P$. The $2$-reachable vertices are double hatched.}
    			\label{fig3}
	\end{figure}

Regions of reachability reflect the nature of pebbling in more detail than values of vertices. We will later analyse the average value of a vertex in a region in a grid. But first, we make a few general observations about regions.

\begin{Obs}
If $P$ is a solvable distribution on a graph $G$ then there exists a non-negative integer $k$ such that $V(G)$ can be partitioned as $R_1 \cup R_2 \cup \ldots R_k \cup S$ where $R_1, \ldots, R_k$ are regions and $S$ is the set of exactly $1$-reachable vertices with a pebble on them.
\end{Obs}

\begin{proof}
We need to show that regions are disjoint and that if $x$ is a vertex not belonging to any region, then $x$ has a pebble on it and is not $2$-reachable.

Assume for contradiction that $v \in R_i \cap R_j$ for $i \neq j$. By maximality, $x$ is not $2$-reachable. That means it has a $2$-reachable neighbour $u_i$ in $R_i$ and a $2$-reachable neighbour $u_j$ in $R_j$. There is no sequence of moves that would result in two pebbles on both $u_i$ and $u_j$ as $x$ is not $2$-reachable. Consider a shortest sequence of moves from $P$ that results in $2$ pebbles on $u_i$ and consider the first move that causes $u_j$ not to be $2$-reachable anymore. Call the vertex from which the move begins $w$. Then there is a path on $2$-reachable vertices between $w$ and $u_i$. There is also a path on $2$-reachable vertices between $w$ and $u_j$. But that means that there is a path on $2$-reachable vertices between $u_i \in R_i$ and $u_j \in R_j$, a contradiction.

On the other hand, let $x$ be a vertex not belonging to any region. Then $x$ is not $2$-reachable. Moreover, it is not a neighbour of any $2$-reachable vertex and therefore no sequence of moves can lead to a pebble being added to $x$. At the same time, $x$ is $1$-reachable (since $P$ is solvable), and therefore $P(x)=1$.
\end{proof}

Next, we shall show that when studying pebbling numbers, we may assume that in a solvable distribution on a connected graph, there is no isolated vertex with exactly one pebble.

\begin{lemma}\label{NoOne}
Let $P$ be a solvable distribution on a connected graph $G$ on at least two vertices. Then there exists a solvable distribution $Q$ on $G$ such that $V(G)$ can be partitioned into regions and $|Q|=|P|$.
\end{lemma}

\begin{proof}
Write $P=R_1 \cup R_2 \cup \ldots R_k \cup S$ where $R_1, \ldots, R_k$ are regions and $S$ is the set of exactly $1$-reachable vertices with a pebble on them. We prove the lemma by induction on $|S|$. If $|S|=0$, we put $Q=P$. Otherwise, take any $x \in S$ and any of its neighbours $u$. As $P$ is solvable, $u$ is reachable. Consider a starting distribution $\hat{P}$ where $\hat{P}(x)=0$, $\hat{P}(u)=P(u)+1$ and $\hat{P}(y)=P(y)$ for all other vertices. Then $\hat{P}$ is still solvable as $x$ has a $2$-reachable neighbour, $|\hat{P}|=|P|$, and $|S|$ has decreased. 
\end{proof}

Finally, we observe that the vertices on the boundary of a region have value at least $\frac{3}{2}$.

\begin{Obs}\label{threehalf}
Let $x$ be a vertex with a neighbour from a region other than the one containing $x$. Then $v(x) \geq \frac{3}{2}$.
\end{Obs}

\begin{proof}
By assumption, $x$ has a $2$-reachable neighbour $u$ from the same region and a $1$-reachable neighbour $x$ from a different region. Similarly to the proof of Observation $\ref{NoOne}$, we show that there exists a sequence of moves after which there are $2$ pebbles on $u$ and a pebble on $x$, therefore $v(x) \geq \frac{3}{2}$. 
\end{proof}

\section{A lower bound on the value of a vertex}

In this section we shall make use of the geometry of the grid $\Lambda_{m,n}$. Let $P$ be any solvable distribution on a grid $\Lambda_{m,n}$, then by the \emph{hemmed pebbling distribution of $P$} we mean the pebbling distribution $P'$ equal to $P$ on the inside of the grid and with two pebbles more than $P$ on every vertex on the boundary. Obviously $|P'| \leq |P|+4(m+n)$, so for the rest of the paper we will focus only on the hemmed pebbling distribution $P'$ and on getting a lower bound on $|P'|$, which will indeed lead to a lower bound on $|P|$. 

\begin{lemma}
Let $P$ be a solvable pebbling distribution on $\Lambda_{m,n}$ and $P'$ be its hemmed pebbling distribution. Then in $P'$ every vertex $X$ satisfies $v(X) \geq 4/3$. \label{FirstLP}
\end{lemma}

\begin{proof}
	
 Assume $X$ is a vertex with minimal value. If $X$ is $2$-reachable, then $v(X) \geq 2$ and we are done. Otherwise $X$ is exactly $1$-reachable. Since all the vertices on the edges of the grid are $2$-reachable, we can assume $X$ is in the inside of the grid. Also, according to Lemma \ref{NoOne}, we can assume that $X$ does not have a pebble on it. Since $X$ is exactly $1$-reachable, then one of its neighbours is $2$-reachable, say $v_2$ by symmetry (see Figure \ref{fig1}). Moreover, since $X$ achieves the minimum value of the grid, then the value of $v_6$ is not smaller than the value of $X$. Hence we have:
 
 \begin{itemize}
    \item $v(v_2)=2A+2B+2C+D/2+E/2+F/2+G/2+H/2 \geq 2,$
    \item $A/2+B/2+C/2+D/2+2E+2F+2G+H/2=v(v_6) \geq v(X) = A+B+C+D+E+F+G+H.$
 \end{itemize}
 
\begin{figure}[htbp]\centering
    			\includegraphics[height=3.5cm]{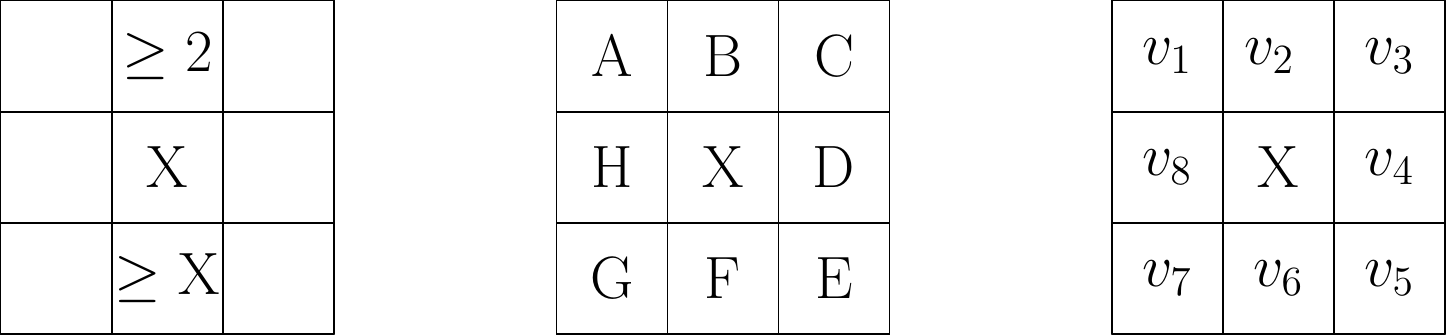}
    			\caption{A visual representation of values, contributions and vertices for proof of Lemma \ref{FirstLP}.}
    			\label{fig1}
	\end{figure}
	
Where, as shown in Fig \ref{fig1}, $A$ is the contribution to the value of $X$ of the top left part of the grid compared to $X$, and $B$, $C$, $D$, $E$, $F$, $G$ and $H$ are defined analogously. \\

Hence, we are looking for the minimum of $A+B+C+D+E+F+G+H$ under the constraints:

\begin{itemize}
    \item $2A+2B+2C+D/2+E/2+F/2+G/2+H/2 \geq 2,$
    \item $-A/2-B/2-C/2-D/2+E+F+G-H/2 \geq 0.$
 \end{itemize}
 
Note that the coefficient of $A$ is not smaller than the coefficient of $B$ in both constraints equations, so for every optimal solution $(A,B,C,D,E,F,G,H)$, $(A+B,0,C,D,E,F,G,H)$ is also an optimal solution. Hence, we can assume $B=0$. Similarly, we can assume $D=F=H=0$. Hence the problem is equivalent to finding the minimum of $A+C+E+G$ under the constraints:

\begin{itemize}
    \item $2A+2C+E/2+G/2 \geq 2,$
    \item $-A/2-C/2+E+G \geq 0.$
 \end{itemize}
 
Taking $S_1=A+C$ and $S_2=E+G$, the problem is equivalent to finding the minimum of $S_1+S_2$ under the constraints:

\begin{itemize}
    \item $2S_1 +\frac{1}{2}S_2 \geq 2,$
    \item $S_2 \geq \frac{1}{2}S_1.$
 \end{itemize}
 
But \[\begin{aligned}
2 \leq 2S_1 +\frac{1}{2}S_2&= \frac{1}{2}S_1 + \frac{3}{2}S_1 +\frac{1}{2}S_2 \\
&\leq S_2 + \frac{3}{2}S_1 +\frac{1}{2}S_2 \\
&= \frac{3}{2}( S_1 +S_2). \\
	\end{aligned}\]\\
	
Consequently, $S_1+S_2 \geq 4/3$, which proves $v(X) \geq 4/3$.
 \end{proof}

\section{A lower bound on the average value in a region}

In this section we analyze regions in grids. First we observe that we have control over the proportion of $2$-reachable vertices in a region (which will allow us to find a better lower bound on the average value of a vertex inside a region). This can be easily proven by induction:

\begin{Obs}\label{boundAlpha}
A region in $\Lambda_{m,n}$ with $k$ $2$-reachable vertices contains at most $3k+2$ vertices.
\end{Obs}

The previous observation claims that the proportion of $2$-reachable vertices is roughly at least $1/3$. Using that and our lower bound on the value of a vertex, we could already get a lower bound on $\pi(\Lambda_{m,n})$, but we can actually improve our bound on the value of the $2$-reachable vertices and get a slightly better lower bound. For this purpose, we will use the following notation: for a $2$-reachable vertex $X$, the \emph{extra value of $X$} is $e(X)=v(X)-2$.

\begin{lemma}\label{AddValue}
Let $P$ be a solvable pebbling distribution on $\Lambda_{m,n}$ and $P'$ be its hemmed pebbling distribution. Let $X$ be a vertex that is $2$-reachable in a region that contains at least two $2$-reachable vertices for the pebbling distribution $P'$. Let $p$ be the number of pebbles on $X$ at the starting configuration. Then:

\begin{enumerate}[(i)]
    \item If $p \geq 3$, then $e(X) \geq p-2.$
    \item If $p=2$, then $e(X) \geq 2/3.$
    \item If $p=1$, then $e(X) \geq 11/75$ and $X$ has a neighbour $Y$ such that $e(Y) \geq 1/2.$
    \item If $p=0$, then either $X$ has a neighbour $Y$ such that $e(Y) \geq 2$ or two neighbours $Y_1$ and $Y_2$ such that $e(Y_1) \geq 1/2$ and $e(Y_2) \geq 1/2.$
\end{enumerate}

\end{lemma}

\begin{proof}
The first point is obvious. For (ii), if $X$ is on the boundary, the result is obvious. Otherwise, we know that there is a neighbour of $X$ that is $2$-reachable, since the region contains at least two $2$-reachable vertices. Without loss of generality assume this neighbour is $v_2$. But $v_5$ and $v_7$ have values at least $4/3$ according to Lemma \ref{FirstLP}. Hence, we are looking for the minimum of $A+B+C+D+E+F+G+H$ under the constraints:

\begin{itemize}
    \item $v(v_2)=2A+2B+2C+D/2+E/2+F/2+G/2+H/2 \geq 2-1,$
    \item $v(v_5)=A/4+B/4+C+D+4E+F+G+H/4 \geq 4/3-1/2,$
    \item $v(v_7)=A+B/4+C/4+D/4+E+F+4G+H \geq 4/3-1/2.$
 \end{itemize}
 
\begin{figure}[htbp]\centering
    			\includegraphics[height=3.5cm]{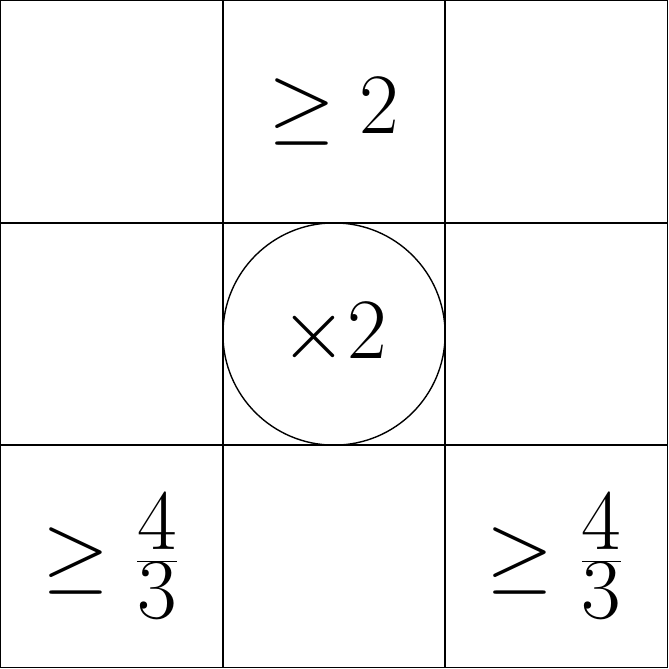}
    			\caption{Proof of Lemma \ref{AddValue}, case $p=2$.}
    			\label{fig2}
	\end{figure}
 
As in the proof of Lemma \ref{FirstLP}, we can suppose $B=D=F=H=0$, which gives the equivalent problem of finding the minimum of $A+C+E+G$ under the constraints:

\begin{enumerate}[(1)]
    \item $2A+2C+E/2+G/2 \geq 1$,
    \item $A/4+C+4E+G \geq 5/6$,
    \item $A+C/4+E+4G \geq 5/6$.
 \end{enumerate}
 
Now $(1)+\frac{2}{5}(2)+\frac{2}{5}(3)$ gives $\frac{5}{2}(A+C+E+G) \geq \frac{5}{3}$, which proves $e(X) \geq 2/3$. \\

For (iii), we know that there exists a neighbour $Y$ of $X$ such that $Y$ is $2$-reachable without using the pebble on $X$. Hence, $e(Y) \geq 1/2$. Moreover, without loss of generality $Y$ is $v_2$, and we know that $v_5$ and $v_7$ have values at least $4/3$ according to Lemma \ref{FirstLP}. Thus, we are looking for the minimum of $A+B+C+D+E+F+G+H$ under the constraints:

\begin{itemize}
    \item $v(v_2)=2A+2B+2C+D/2+E/2+F/2+G/2+H/2+1/2 \geq 2+1/2,$
    \item $v(v_5)=A/4+B/4+C+D+4E+F+G+H/4 \geq 4/3-1/4,$
    \item $v(v_7)=A+B/4+C/4+D/4+E+F+4G+H \geq 4/3-1/4.$
 \end{itemize}
 
As in the proof of Lemma \ref{FirstLP}, we can suppose $B=D=F=H=0$, which gives the equivalent problem of finding the minimum of $A+C+E+G$ under the constraints:

\begin{enumerate}[(1)]
    \item $2A+2C+E/2+G/2 \geq 2$,
    \item $A/4+C+4E+G \geq 13/12$,
    \item $A+C/4+E+4G \geq 13/12$.
 \end{enumerate}
 
Now $(1)+\frac{2}{5}(2)+\frac{2}{5}(3)$ gives $\frac{5}{2}(A+C+E+G) \geq \frac{43}{15}$. Hence $A+C+E+G \geq \frac{86}{75} $ which proves $e(X) \geq \frac{11}{75}$. \\

For (iv), since $X$ is $2$-reachable and has no pebble on it at the starting configuration, it either has a neighbour $Y$ that is $4$-reachable, or two neighbours $Y_1$ and $Y_2$ that are $2$-reachable at the same time. In the first case, $e(Y) \geq 2$. In the second case, since $Y_1$ and $Y_2$ are at distance $2$ of each other, they will give to each other a contribution of $1/2$. Hence $e(Y_1) \geq 1/2$ and $e(Y_2) \geq 1/2$.
\end{proof}

\begin{lemma}\label{LowerBoundRegion}
Let $P$ be a solvable pebbling distribution on $\Lambda_{m,n}$ and $P'$ be its hemmed pebbling distribution. In $P'$, let $R$ be a region, $k \geq 2$ the number of vertices that are $2$-reachable in $R$ and $N$ the total number of vertices in $R$. Then the average value of $R$, denoted by $A(R)$, satisfies: 

$$A(R) \geq \frac{(2+\frac{50}{353})k + \frac{4}{3}(N-k)}{N}.$$

Moreover, if the grid contains more than $2$ regions:

$$A(R) \geq \frac{(2+\frac{50}{353})k + \frac{4}{3}(N-k) +\frac{2}{3}}{N}.$$

\end{lemma}

\begin{proof}
 We apply Lemma \ref{AddValue} to each $2$-reachable vertex. We cannot add the extra values given in the cases $p=0$, $p=1$ and $p \geq 2$ because in some cases, we add an extra value to the vertex $X$, and in some other cases, we add an extra value to a neighbouring vertex, so we have to be careful not to double-count some extra values. Let $x$, $y$ and $z$ be respectively the number of $2$-reachable vertices with $p=0$, $p=1$ and $p \geq 2$. Considering the extra values on the vertices themselves, we get an overall extra value of at least $e_1=\frac{2}{3}z+\frac{11}{75}y$.

Considering the extra values of the neighbours of $X$ in the case $p=0$, since each neighbour could be counted $4$ times, we get an overall extra value of at least $e_2=\frac{1}{4}x$.

Considering the extra values of the neighbours in the case $p=1$, by taking $Y$ to be a neighbour of $X$ that can get $2$ pebbles in the minimum number of moves, we make sure that we do not double count extra values, since in the case $t > 1$ vertices $X_1, \dots, X_t$  share the same $Y$, there is a sequence of moves that put 2 pebbles on $Y$ without using the pebble on $X_1, \dots, X_t$, so $e(Y) \geq t/2$. Hence we get an overall extra value of at least $e_3=\frac{1}{2}y$.

We have $x+y+z=k$, and we can always choose to add the biggest quantity between $e_1$, $e_2$ and $e_3$, so the worst-case scenario will be when those $3$ quantities are equal. When this happens, we have $x=2y$ and $z=\frac{53}{100}y$, so $2y+y+\frac{53}{100}y=k$, which gives $y=\frac{100}{353}k$. Hence the overall extra value is always at least $\frac{50}{353}k$, which proves that:

$$A(R) \geq \frac{(2+\frac{50}{353})k + \frac{4}{3}(N-k)}{N}.$$

Moreover, in the case where the grid is split in at least $2$ regions, there exists at least $4$ vertices in $R$ that have at least one neighbour from a different region. Those vertices have values at least $3/2$ according to Observation \ref{threehalf}. Hence, we get:

$$A(R) \geq \frac{(2+\frac{50}{353})k + \frac{4}{3}(N-k) +\frac{2}{3}}{N},$$

as claimed.
\end{proof}

\begin{lemma}\label{Casediffregions}
Let $P$ be a solvable pebbling distribution on $\Lambda_{m,n}$ and $P'$ be its hemmed pebbling distribution. Then in $P'$, let $R$ be a region, and suppose the grid is split in at least $2$ regions. Then $$A(R) \geq \frac{5092}{3177}.$$
\end{lemma}

\begin{proof}

Let $k$ be the number of $2$-reachable vertices inside $R$ and $N$ be the total number of vertices inside $R$. If $k \geq 2$, then $N \leq 3k+2$ as stated in Observation \ref{boundAlpha}, and the bound in Lemma \ref{LowerBoundRegion} gives: 

\[\begin{aligned}
A(R) &\geq \frac{k (2+\frac{50}{353}) + (2k+2) \frac{4}{3} + \frac{2}{3}}{3k+2} = \frac{5092}{3177}+\frac{406}{3177(3k+2)} \geq \frac{5092}{3177}. \\
	\end{aligned}\]\\
	
If $k=1$, then let $X$ be the $2$-reachable vertex and use the classic notation: $v_2$, $v_4$, $v_6$ and $v_8$ have values at least $3/2$ according to \ref{threehalf}, and since $v_1$ is reachable within its own region, it will give at least $1/4$ extra contribution to $X$. Hence, $A(R) \geq \frac{4*3/2+2+1/4}{5} \geq \frac{5092}{3177}$. 

\end{proof}

\begin{theorem}
Let $P$ be a solvable pebbling distribution on $\Lambda_{m,n}$, then: 

$$|P| \geq \frac{5092}{28593}nm +O(n+m).$$
\end{theorem}

\begin{proof}
Recall that $P'$ is the hemmed pebbling distribution where we have added $2$ pebbles on every vertex that lies on the edge. Then according to Lemma \ref{LowerBoundRegion}, we have: 

$$9|P| \geq 9|P'|+ O(n+m) \geq \sum_{R} A(R)|R| + O(n+m).$$

If the grid is split in at least $2$ regions, then according to Lemma \ref{Casediffregions}, we have: 
$$\sum_{R} A(R)|R| \geq \sum_{R} \frac{5092}{3177}|R|=\frac{5092}{3177}nm.$$

Hence $|P| \geq \frac{5092}{28593}nm +O(n+m)$ as wanted. \\

If the grid contains only one region $R$, then using Lemma \ref{LowerBoundRegion}, we get: 

\[\begin{aligned}
A(R) &\geq \frac{(2+\frac{50}{353})k + \frac{4}{3}(N-k)}{N} \\
&\geq \frac{5092}{3177} + O(1/N) \\
&= \frac{5092}{3177} + O(n+m). \\
	\end{aligned}\]
	
Which finally gives:

$$|P| \geq \frac{5092}{28593}nm +O(n+m),$$

as claimed.
 \end{proof}

\section{Concluding remarks and open problems}

As was mentioned in the introduction, the new lower bound on $\pi(\Lambda_{m,n})$ is far from the best known upper bound. One of the reasons is that in our proof, we only distinguish between $2$-reachable and exactly $1$-reachable vertices. A possible improvement could come from taking into account $3$-reachable, or even $4$-reachable vertices, as the best known construction (\cite{UpperBoundGrid}) contains linearly many (in terms of the total number of vertices) vertices with $4$ pebbles on them.

Another line of improvement would be to improve the lower bound of the value of a vertex from $\frac{4}{3}$ to something higher. In view of Observation \ref{threehalf}, the authors believe the following should be true:

\begin{conjecture}
Let $P$ be a solvable pebbling distribution on the grid. Then for all vertices $X$ that do not lie on the boundary of the grid, $v(X) \geq \frac{3}{2}$.
\end{conjecture}

The reader would have probably recognised Linear Programming equations in Lemma \ref{FirstLP} and Lemma \ref{AddValue}. It may seem that some equations could be optimised by adding more constraints on the values of remaining neighbouring vertices. However, adding them will not improve the lower bound \textemdash the authors' approach was to begin with all the constraints on the $3 \times 3$ square centered in $X$ and then use a computer software to delete useless ones. What remained was a small number of constraints that was possible to solve 'by hand'. \\

We finish the paper by mentioning two variations of pebbling a square grid for which we could achieve lower bounds using the arguments presented in this paper. First, we could consider pebbling of a $k$-dimensional grid.

A different variation is to consider \emph{$k$-pebbling moves}: removing $k$ pebbles from a vertex to add a pebble on a neighbouring vertex. What we called a pebbling move is then a $2$-pebbling move. Once again, given a graph $G$ and a pebbling distribution $P$, we will call a vertex $y$ reachable if there is a sequence of $k$-pebbling moves starting in $P$ and ending in a distribution $Q$ with $Q(y) \geq 1$. A pebbling distribution $P$ on a graph $G$ is \emph{$k$-solvable} if every vertex is reachable. For a graph $G$, we define the \emph{optimal k-pebbling number $\pi_{k}(G)$ of $G$} as the minimal total number of pebbles among all $k$-solvable distributions on $G$. For any connected graph $G$, $\pi_{1}(G)=1$. For a grid $\Lambda_{m,n}$, the situation is also simple when $k\geq 5$ in view of the following lemma:

\begin{lemma}\label{smoothening}
Let $k \in \mathbb{N}$ and $P$ be a $k$-solvable distribution on $\Lambda_{m,n}$ such that there is a vertex $x$ with $P(x) \geq k+1$. Then $\hat{P}$, where
\begin{itemize}
    \item $\hat{P}(x)=P(x)-k$,
    \item $\hat{P}(z)=P(z)+1$ if $z$ is a neighbour of $x$, and
    \item $\hat{P}(z)=P(z)$ if $d(x,z)\geq 2$,
\end{itemize}
is also a $k$-solvable distribution.
\end{lemma}

\begin{proof}
Let $y$ be any vertex of $\Lambda_{m,n}$. We show it is reachable in $\hat{P}$. This is clearly true if $d(x,y) \leq 1$.

Suppose $y$ is such that $d(x,y) \geq 2$. Fix a sequence $\Sigma$ of moves that begins in distribution $P$ and results in a distribution with a pebble on $y$. If none of the moves begins in $x$, the same sequence of moves witnesses the reachability of $y$ in $\hat{P}$.

Otherwise, consider the first move $M$ in $\Sigma$ that removes pebbles from $x$. Let $Q$ be the distribution after all the moves from $\Sigma$ up to $M$ (inclusive) have been made. We claim that $\Sigma$ without $M$ is a sequence of moves witnessing the reachability of $y$ in $\hat{P}$.

Indeed, as $\hat{P}(z) \geq {P}(z)$ for all $z \neq x$, all moves up to $M$ (exclusive) are possible also when beginning in $\hat{P}$. Call the distribution after all these moves have been made $\hat{Q}$. Then $\hat{Q}(z) \geq Q(z)$ for all vertices $z$ (with a sharp inequality for all neighbours of $x$ but one). Therefore, all the remaining moves from $\Sigma$ without $M$ can also be made, after which $y$ has a pebble on it. 
\end{proof}

\begin{Prop}
For $k \geq 5$, $\pi_k(\Lambda_{m,n}) = nm$.
\end{Prop}

\begin{proof}
The pebbling distribution with $1$ pebble on each vertex is $k$-solvable, thus $\pi_k(\Lambda_{m,n}) \leq nm$.

For the other direction, consider a $k$-solvable distribution $P$. We show that there is a $k$-solvable distribution $Q$ with at least $1$ pebble on each vertex and satisfying $|Q|\leq |P|$. Note that this implies $\pi_{k}(\Lambda_{m,n})\geq nm$, finishing the proof.

Let $x$ be a vertex such that $P(x)=0$. Since $x$ is reachable, there exists a sequence $\Sigma$ of $l$ $k$-pebbling moves that results in a pebble on $x$. For each of the $l$ moves, starting with the first move, we use Lemma \ref{smoothening} to get a $k$-solvable distribution $\hat{P}$. Since each vertex has at most $4$ neighbours, we have $|\hat{P}|\leq |P|$. Note that the number of unoccupied vertices in $\hat{P}$ is not smaller than in $P$. After $l$ uses of Lemma \ref{smoothening}, we arrive at a $k$-solvable distribution with at most $|P|$ pebbles and with fewer unoccupied vertices than $P$ has (as $x$ has a pebble on it).

By repeating the process from the previous paragraph we eventually obtain the desired configuration $Q$.
\end{proof}

What would the optimal $k$-pebbling number of $\Lambda_{m,n}$ be for $k \in \{3,4\}$?

\section*{Acknowledgement}

The authors would like to thank Béla Bollobás for suggesting this problem and for his valuable comments.

\bibliographystyle{abbrvnat}  
	\renewcommand{\bibname}{Bibliography}
	\bibliography{bibliography}

\end{document}